\documentclass[12pt,leqno]{article}
\usepackage{amsfonts,amssymb,amsmath,amsthm}
\newtheorem{theorem}{{\sc Theorem}}

\newtheorem{lemma}{Lemma}
\theoremstyle{definition}

\def\RR{\mathbb{R}}
\def\QQ{\mathbb{Q}}
\def\ZZ{\mathbb{Z}}
\def\NN{\mathbb{N}}

\newcommand{\ie}{{\it i.e.}\/ }

\renewcommand{\le}{\leqslant}
\renewcommand{\ge}{\geqslant}
\newcommand{\vv}[1]{{\mathbf{#1}}}

\newcommand{\dist}{\operatorname{dist}}

\newcommand{\Sprindzuk}{Sprind\v{z}uk }
\newcommand{\w}{W_n(\Psi)}
\newcommand{\wirr}{W^*_n(\Psi)}

\newcommand{\wb}{W_{{\rm big}}(\Psi)}
\newcommand{\wm}{W_{{\rm med}}(\Psi)}
\newcommand{\ws}{W_{{\rm small}}(\Psi)}
\newcounter{mysection}
\newcommand{\mysection}[1]{\refstepcounter{mysection}\par\bigskip{\bf\arabic{mysection}.\,\,#1.} }

\begin{document}

\begin{center}
{\large\bf On a theorem of V.~Bernik in the metrical theory of
Diophantine approximation}\\[1ex]
by\\[1ex]  {\sc V.~Beresnevich~(Minsk)}%
\footnote{The work has been supported by EPSRC grant
GR/R90727/01\par {\ \ \it Key words and phrases}\/: Diophantine
approximation, Metric theory of Diophantine approximation, the
problem of Mahler\par {\ \ 2000 \it Mathematics Subject
Classification}\/: 11J13, 11J83, 11K60}
\end{center}

\mysection{Introduction}
We begin by introducing some notation: $\#S$ will denote the
number of elements in a finite set $S$; the Lebesgue measure of a
measurable set $S\subset\RR$ will be denoted by $|S|$; $P_n$ will
be the set of integral polynomials of degree $\le n$. Given a
polynomial $P$, $H(P)$ will denote the height of $P$, \ie the
maximum of the absolute values of its coefficients; $P_n(H)=\{P\in
P_n:H(P)=H\}$. The symbol of Vinogradov $\ll$ in the expression
$A\ll B$ means $A\le CB$, where $C$ is a constant. The symbol
$\asymp$ means both $\ll$ and $\gg$. Given a point $x\in\RR$ and a
set $S\subset\RR$, $\dist(x,S)=\inf\{|x-s|:x\in S\}$. Throughout,
$\Psi$ will be a positive function.

\smallskip

{\bf Mahler's problem.}
In 1932 K.~Mahler \cite{Mahler-1932a} introduced a classification
of real numbers $x$ into the so-called classes of $A,S,T$ and $U$
numbers according to the behavior of $w_n(x)$ defined as the
supremum of $w>0$ for which
$$
|P(x)|<H(P)^{-w}
$$
holds for infinitely many $P\in P_n$. By Minkovski's theorem on
linear forms, one readily shows that $w_n(x)\ge n$ for all
$x\in\RR$. Mahler \cite{Mahler-1932b} proved that for almost all
$x\in\RR$ (in the sense of Lebesgue measure) $w_n(x)\le 4n$, thus
almost all $x\in\RR$ are in the $S$-class. Mahler has also
conjectured that for almost all $x\in\RR$ one has the equality
$w_n(x)=n$. For about 30 years the progress in Mahler's problem
was limited to $n=2$ and $3$ and to partial results for $n>3$.
V.~Sprindzuk proved Mahler's conjecture in full (see
\cite{Sprindzuk-1969}).

\smallskip {\bf A.~Baker's conjecture.}
Let $\w$ be the set of $x\in \RR$ such that there are infinitely
many $P\in P_n$ satisfying
\begin{equation}\label{e:001}
|P(x)|<\Psi(H(P)).
\end{equation}
A.~Baker \cite{Baker-1966} has improved Sprind\v{z}uk's theorem by
showing that
$$
|\w|=0\text{ if }\sum_{h=1}^\infty\Psi^{1/n}(h)<\infty\text{ and
$\Psi$ is monotonic.}
$$
He has also conjectured a stronger statement proved by V.~Bernik
\cite{Bernik-1989} that $|\w|=0$ if the sum
\begin{equation}\label{e:002}
  \sum_{h=1}^\infty h^{n-1}\Psi(h)
\end{equation}
converges and $\Psi$ is monotonic. Later V.~Beresnevich
\cite{Beresnevich-1999c} has shown that $|\RR\setminus\w|=0$ if
(\ref{e:002}) diverges and $\Psi$ is monotonic. We prove

\begin{theorem}\label{theorem1}
Let $\Psi:\RR\to\RR^+$ be arbitrary function $($not necessarily
monotonic$)$ such that the sum\/ {\rm (\ref{e:002})}\/ converges.
Then $|\w|=0$.
\end{theorem}

Theorem~\ref{theorem1} is no longer improvable as, by
\cite{Beresnevich-1999c}, the convergence of (\ref{e:002}) is
crucial. Notice that for $n=1$ the theorem is simple and known
(see, for example, \cite[p.\,121]{Cassels1}). Therefore, from now
on we assume that $n\ge 2$.

\mysection{Subcases of Theorem~\ref{theorem1}}
Let $\delta>0$. We define the following 3 sets denoted by $\wb$,
$\wm$ and $\ws$ consisting of $x\in\RR$ such that there are
infinitely many $P\in P_n$ simultaneously satisfying (\ref{e:001})
and one of the following inequalities
\begin{equation}\label{e:003}
    1\le|P'(x)|,
\end{equation}
\begin{equation}\label{e:004}
    H(P)^{-\delta}\le |P'(x)|<1,
\end{equation}
\begin{equation}\label{e:005}
    |P'(x)|< H(P)^{-\delta}
\end{equation}
respectively. Obviously $\w=\wb\cup\wm\cup\ws$. Hence to prove
Theorem~\ref{theorem1} it suffices to show that each of the sets
has zero measure.

Since sum (\ref{e:002}) converges,
$H^{n-1}\Psi(H)$ tends to $0$ as $H\to\infty$. Therefore,
\begin{equation}\label{e:006}
  \Psi(H)=o(H^{-n+1})\text{ as }H\to\infty.
\end{equation}

\mysection{The case of a big derivative}
The aim of this section is to prove that $|\wb|=0$. Let $B_{n}(H)$
be the set of $x\in \RR$ such that there exists a polynomial $P\in
P_n(H)$ satisfying (\ref{e:003}). Then
\begin{equation}\label{e:007}
\textstyle\wb=\bigcap\limits_{N=1}^\infty\
\bigcup\limits_{H=N}^\infty B_n(H).
\end{equation}
Now $|\wb|=0$ if $|\wb\cap I|=0$ for any open interval $I\subset
\RR$ satisfying
\begin{equation}\label{e:008}
  0<c_0(I)=\inf\{|x|:x\in I\} < \sup\{|x|:x\in I\}=c_1(I)<\infty.
\end{equation}
Therefore we can fix an interval $I$ satisfying (\ref{e:008}).

By (\ref{e:007}) and the Borel-Cantelli Lemma, $|\wb\cap I|=0$
whenever
\begin{equation}\label{e:009}
\sum_{H=1}^\infty |B_n(H)\cap I|<\infty.
\end{equation}

By the convergence of (\ref{e:002}), condition (\ref{e:009}) will
follow on showing that
\begin{equation}\label{e:010}
|B_n(H)\cap I|\ll H^{n-1}\Psi(H)
\end{equation}
with the implicit constant in (\ref{e:010}) independent of $H$.

Given a $P\in P_n(H)$, let $\sigma(P)$ be the set of
$x\in I$ satisfying (\ref{e:003}). Then
\begin{equation}\label{e:011}
B_{n}(H)\cap I=\bigcup\nolimits_{P\in P_n(H)}\sigma(P).
\end{equation}

\begin{lemma}\label{lem2}
Let $I$ be an interval with endpoints $a$ and $b$. Define the
following sets
$
I''=[a,a+4\Psi(H)]\cup[b-4\Psi(H),b]\mbox{ \ \ and \ \ }
I'=I\setminus I''.
$
Then for all sufficiently large $H$ for any $P\in P_n(H)$ such
that $\sigma(P)\cap I'\not=\emptyset$, for any
$x_0\in\sigma(P)\cap I'$ there exists $\alpha\in I$ such that
$P(\alpha)=0$, $|P'(\alpha)|>|P'(x_0)|/2$ and
$|x_0-\alpha|<2\Psi(H)\,|P'(\alpha)|^{-1}$.
\end{lemma}

The proof of this Lemma nearly coincides with the one of Lemma~1
in \cite{Beresnevich-1999c} and is left for the reader. There will
be some changes to constants and notation and one also will have
to use (\ref{e:006}).


Given a polynomial $P\in P_n(H)$ and a real number $\alpha$ such
that $P'(\alpha)\not=0$, define $\sigma(P;\alpha)=\left\{x\in
I:|x-\alpha|<2\Psi(H)|P'(\alpha)|^{-1}\right\}.$ Let $I'$ and
$I''$ be defined as in Lemma~\ref{lem2}. For every polynomial
$P\in P_n(H)$, we define the set
$$
Z_I(P)=\{\alpha\in I:P(\alpha)=0\mbox{ and }|P'(\alpha)|\ge1/2\}.
$$
By Lemma~\ref{lem2}, for any $P\in P_n(H)$ we have the inclusion
\begin{equation}\label{e:012}
\sigma(P)\cap I'\subset\bigcup\nolimits_{\alpha\in
Z_I(P)}\sigma(P;\alpha).
\end{equation}

Given $k\in\ZZ$ with $0\le k\le n$, define
$$
P_n(H,k)=\{P=a_nx^n+\dots+a_0\in P_n(H): a_k=0\}
$$
and for $R\in P_n(H,k)$ let
$
P_n(H,k,R)=\{P\in P_n(H):P-R=a_kx^k\}.
$
It is easily observed that
\begin{equation}\label{e:013}
P_n(H)=\bigcup_{k=0}^n \ \bigcup\nolimits_{R\in P_n(H,k)}
P_n(H,k,R)
\end{equation}
and
\begin{equation}\label{e:014}
\#P_n(H,k)\ll  H^{n-1}\text{ \ \ for every }k.
\end{equation}
Taking into account (\ref{e:011}), (\ref{e:013}), (\ref{e:014}) and that
$|I''|\ll\Psi(H)$, it now becomes clear that to prove
(\ref{e:010}) it is sufficient to show that for every fixed $k$
and fixed $R\in P_n(H,k)$
\begin{equation}\label{e:015}
\Big|\bigcup_{P\in P_n(H,k,R)}\sigma(P)\cap I'\,\Big|\ll \Psi(H).
\end{equation}

Let $k$ and $R$ be fixed. Define the rational function $\tilde
R(x)=x^{-k}R(x)$. By (\ref{e:008}), there exists a collection of
intervals $[w_{i-1},w_i)\subset I$ $(i=1,\dots,s)$, which do not
intersect pairwise and cover $I$, such that $\tilde R(x)'$ is
monotonic and does not change the sign on every interval
$[w_{i-1},w_i)$. It is clear that $s$ depends on $n$ only. Let
$Z_{I,R}=\bigcup\nolimits_{P\in P_n(H,k,R)}Z_I(P)$,
$k_i=\#(Z_{I,R}\cap[w_{i-1},w_i))$ and $Z_{I,R}\cap
[w_{i-1},w_i)=\{\alpha_i^{(1)},\dots,\alpha_i^{(k_i)}\}$, where
$\alpha^{(j)}_i<\alpha^{(j+1)}_i$. Given a $P\in P_n(H,k,R)$, we
obviously have the identity
$$
\frac{x^kP'(x)-kx^{k-1}P(x)}{x^{2k}}=\left(\frac{P(x)}{x^k}\right)'=
\tilde R(x)'.
$$
Taking $x$ to be $\alpha\in Z_I(P)$ leads to
$\frac{P'(\alpha)}{\alpha^{k}} = \tilde R(\alpha)'$. By
(\ref{e:008}), $|P'(\alpha)|\asymp|\tilde R(\alpha)'|$. Now, by
Lemma~\ref{lem2}, $|\sigma(P;\alpha)|\ll
\Psi(H)\,|P'(\alpha)|^{-1}\ll \Psi(H)\,|\tilde R'(\alpha)|^{-1}$.

Using (\ref{e:012}), we get
$$
\Big|\bigcup_{P\in P_n(H,k,R)}\sigma(P)\cap I'\Big|\ll \Psi(H)
\sum_{i=1}^s\sum_{j=1}^{k_i}\frac{1}{|\tilde R(\alpha_i^{(j)})|}
$$
Now to show (\ref{e:015}) it suffices to prove that for every $i$
$(1\le i\le s)$
\begin{equation}\label{e:016}
\sum_{j=1}^{k_i}|\tilde R'(\alpha_i^{(j)})|^{-1}\ll 1.
\end{equation}

Fix an index $i$ $(1\le i\le s)$. If $k_i\ge2$ then we can
consider two sequential roots $\alpha_i^{(j)}$ and
$\alpha_i^{(j+1)}$ of two rational functions $\tilde R+a_k^{i,j}$
and $\tilde R+a_k^{i,j+1}$ respectively. For convenience let us
assume that $\tilde R'$ is increasing and positive on
$[w_{i-1},w_i)$. Then $\tilde R$ is strictly monotonic on
$[w_{i-1},w_i)$, and we have $a_k^{i,j}\not=a_k^{i,j+1}$. It
follows that $|a_k^{i,j}-a_k^{i,j+1}|\ge1$. Using The Mean Value
Theorem and the monotonicity of $\tilde R'$, we get
$$
1\le |a_0^{i,j}-a_0^{i,j+1}|=|\tilde R'(\alpha_i^{(j)})-\tilde
R'(\alpha_i^{(j+1)})|=
$$
$$
=|\tilde R'(\tilde{\alpha}_i^{(j)})|\cdot
|\alpha_i^{(j)}-\alpha_i^{(j+1)}| \le|\tilde
R'(\alpha_i^{(j+1)})|\cdot|\alpha_i^{(j)}-\alpha_i^{(j+1)}|,
$$
where $\tilde{\alpha}_i^{(j)}$ is a point between $\alpha_i^{(j)}$
and $\alpha_i^{(j+1)}$. This implies $|\tilde
R'(\alpha_i^{(j+1)})|^{-1}\le\alpha_i^{(j+1)}-\alpha_i^{(j)},$
whence we readily get
$$
\sum_{j=1}^{k_i-1}|\tilde R'(\alpha_i^{(j+1)})|^{-1}\le
\sum_{j=1}^{k_i-1}\left(\alpha_i^{(j+1)}-\alpha_i^{(j)}\right)=
\alpha_i^{(k_i)}-\alpha_i^{(1)}\le w_{i}-w_{i-1}.
$$
The last inequality and $|\tilde R'(\alpha_i^{(1)})|\asymp|
P'(\alpha_i^{(1)})|\gg 1$ yield (\ref{e:016}). It is easily
verified that (\ref{e:016}) holds for every $i$ with $k_i\ge 2$
and is certainly true when $k_i=1$ or $k_i=0$. This completes the
proof of the case of a big derivative.

\mysection{The case of a medium derivative}
As above we fix an interval $I$ satisfying (\ref{e:008}). The
statement $|\wm|=0$ will now follow from $|\wm\cap I|=0$. We will
use the following

\begin{lemma}[see Lemma 2 in \cite{Beresnevich-2002a}]\label{lem1}
Let
$\alpha_0,\dots,\alpha_{k-1},\beta_1,\dots,\beta_k\in\RR\bigcup
\{+\infty\}$ be such that $\alpha_0>0$, $\alpha_j>\beta_j\ge 0$
for $j=1,\dots,k-1$ and\/ $0<\beta_k<+\infty$. Let $f:(a,b)\to\RR$
be a $C^{(k)}$ function such that $
\inf_{x\in(a,b)}|f^{(k)}(x)|\ge \beta_k. $ Then, the set of
$x\in(a,b)$ satisfying
$$
\left\{
\begin{array}{l}
\hspace*{8ex}|f(x)|\le \alpha_0,\\[0.5ex]
\ \beta_j \ \le \ |f^{(j)}(x)|\le \alpha_j\ \ (j=1,\dots,k-1)
\end{array}
\right.
$$
is a union of at most $k(k+1)/2+1$ intervals with lengths at most
$ \min_{0\le i<j\le
k}3^{(j-i+1)/2}(\alpha_i\left/\beta_j\right.)^{1/(j-i)}$. Here we
adopt $\frac{c}{0}=+\infty$ for $c>0$.
\end{lemma}

Given a polynomial $P\in P_n(H)$, we redefine $\sigma(P)$ to be
the set of solutions of (\ref{e:004}). Since $P^{(n)}(x)=n!a_n$,
we can apply Lemma~\ref{lem1} to $P$ with $k=n$ and
$$
\alpha_0=\Psi(H),\ \alpha_1=1,\
\beta_1=\inf_{x\in\sigma(P)}|P'(x)|\ge H^{-\delta},\ \beta_n=1,
$$
$$
\alpha_2=\dots=\alpha_{n-1}=+\infty,\ \beta_2=\dots=\beta_{n-1}=0.
$$
Then we conclude that $\sigma(P)$ is a union of at most
$n(n+1)/2+1$ intervals of length $\ll \alpha_0/\beta_1$. There is
no loss of generality in assuming that the sets $\sigma(P)$ are
intervals as, otherwise, we would treat the intervals of
$\sigma(P)$ separately. We also can ignore those $P$ for which
$\sigma(P)$ is empty. For every $P$ we also define a point
$\gamma_P\in\sigma(P)$ such that $\inf_{x\in\sigma(P)}|P'(x)|\ge
\frac12|P'(\gamma_P)|$. The existence is easily seen. Now we have
\begin{equation}\label{e:017}
|\sigma(P)|\ll {\Psi(H)}{|P'(\gamma_P)|}^{-1}.
\end{equation}
It also follows from the choice of $\gamma_P$ that
\begin{equation}\label{e:018}
H(P)^{-\delta}\le|P'(\gamma_P)|<1.
\end{equation}

Now define expansions of $\sigma(P)$ as follows:
$$
\sigma_1(P):=\{x\in I:\dist(x,\sigma(P))<(H|P'(\gamma_P)|)^{-1}\},
$$
$$
\sigma_2(P):=\{x\in I:\dist(x,\sigma(P))<H^{-1+2\delta}\}.
$$
By (\ref{e:004}), $\sigma_1(P)\subset\sigma_2(P)$. Moreover, it is
easy to see that
\begin{equation}\label{e:019}
\sigma_1(P)\subset\sigma_2(Q)\text{ \ \ for any $Q\in P_n(H)$ with
$\sigma_1(Q)\cap\sigma_1(P)\not=\emptyset$.}
\end{equation}

It is readily verified that
$|\sigma_1(P)|\asymp(H|P'(\gamma_P)|)^{-1}$, and therefore, by
(\ref{e:017}),
$$
|\sigma(P)|\ll |\sigma_1(P)|\, H \Psi(H).
$$

Take any $x\in\sigma_2(P)$. Using the Mean Value Theorem,
(\ref{e:018}) and $|x-\gamma_P|\ll H^{-1+2\delta}$,
we get $ |P'(x)|\le|P'(\gamma_P)|+|P''(\tilde x)((x-\gamma_P)|\ll
1+H\cdot H^{-1+2\delta}\ll H^{2\delta}, $ where $\tilde x$ is
between $x$ and $\gamma_P$. Similarly we estimate $|P(x)|$
resulting in
\begin{equation}\label{e:020}
|P(x)|\ll H^{-1+4\delta},\ \  |P'(x)|\ll H^{2\delta}\text{ \ for
any $x\in\sigma_2(P)$}.
\end{equation}

Now for every pair $(k,m)$ of integers with $0\le k<m\le n$ we
define
$$
P_n(H,k,m)=\{R=a_nx^n+\dots+a_0\in P_n(H):a_k=a_m=0\}
$$
and for a given polynomial $R\in P_n(H,k,m)$ we define
$$
P_n(H,k,m,R)=\{P=R+a_mx^m+a_kx^k\in P_n(H)\}.
$$

The intervals $\sigma(P)$ will be divided into 2 classes of
essential and non-essential intervals. The interval $\sigma(P)$
will be essential if for any choice of $(k,m,R)$ such that $P\in
P_n(H,k,m,R)$ for any $Q\in P_n(H,k,m,R)$ other than $P$ we have $
\sigma_1(P)\cap\sigma_1(Q)=\emptyset. $ For fixed $k$, $m$ and $R$
summing the measures of essential intervals gives
$$
\sum|\sigma(P)|\le H\Psi(H)\sum|\sigma_1(P)|\le H\Psi(H)|I|\ll
H\Psi(H).
$$

As $\#P_n(H,k,m)\ll H^{n-2}$ and there are only $n(n+1)/2$
different pairs $(k,m)$ we obtain the following estimate
$$ \sum_{\text{essential intervals $\sigma(P)$ with }P\in
P_n(H)}|\sigma(P)|\ll H^{n-1}\Psi(H).
$$
Thus, by the Borel-Cantelli Lemma and the convergence of
(\ref{e:002}), the set of points $x$ of $\wm\cap I$ which belong
to infinitely many essential intervals is of measure zero.

Now let $\sigma(P)$ be non-essential. Then, by definition and
(\ref{e:019}) there is a choice of $k,m,R$ such that $P\in
P_n(H,k,m,R)$ and there is a $Q\in P_n(H,k,m,R)$ different from
$P$ such that
$$
\sigma(P)\subset\sigma_1(P)\subset\sigma_2(P)\cap\sigma_2(Q).
$$
On the set $\sigma_2(P)\cap\sigma_2(Q)$ both $P$ and $Q$ satisfy
(\ref{e:020}) and so does the difference
$P(x)-Q(x)=b_mx^m+b_kx^k$. It is not difficult to see that
$b_m\not=0$ if $H$ is big enough. Therefore using (\ref{e:020})we
get
\begin{equation}\label{e:021}
\Big|x^{m-k}+\frac{b_k}{b_m}\Big|\ll
\frac{H^{-1+4\delta}}{|b_m|}\le H^{-1+4\delta}\quad\text{and}\quad
\max\{|b_m|,|b_k|\}\ll H^{2\delta}.
\end{equation}
Now let $x$ belongs to infinitely many non-essential intervals.
Without loss of generality we assume that $x$ is transcendental as
otherwise it belongs to a countable set, which is of measure zero.
Therefore (\ref{e:021}) is satisfied for infinitely many
$b_m,b_k\in\ZZ$. Hence, the inequality
$$
\Big|x^{m-k}-\frac pq\Big|< q^{-\frac{1-5\delta}{2\delta}}
$$
holds for infinitely many $p,q\in\ZZ$. Taking $\delta=\frac1{10}$
so that $\frac{1-5\delta}{2\delta}$ becomes $2+\delta$ and
applying standard Borel-Cantelli arguments (see
\cite[p.\,121]{Cassels1}) we complete the proof of the case of a
medium derivative for non-essential intervals.

\mysection{The case of a small derivative}
In this section we prove that $|\ws|=0$. We will make use of
Theorem 1.4 in \cite{BernikKleinbockMargulis-2001a}. By taking
$d=1$, $\vv f=(x,x^2,\dots,x^n)$, $U=\RR$, $T_1=\ldots=T_n=H$,
$\theta=H^{-n+1}$, $K=H^{-\delta}$ in that theorem, we arrive at

\begin{theorem}\label{theorem2}
Let $x_0\in \RR$ and
$\delta'=\frac{\min(\delta,n-1)}{(n+1)(2n-1)}$. Then there exists
a finite interval $I_0\subset \RR$ containing $x_0$ and a constant
$E > 0$ such that
$$
\Big|\bigcup_{P\in P_n,\,0<H(P)\le H}\Big\{x\in I_0: |P(x)|  <
H^{-n+1}, |P'(x)|< H^{-\delta} \Big\}\Big|\le E H^{-\delta'}.
$$
\end{theorem}

In particular Theorem~\ref{theorem2} implies that for any
$\delta>0$ the set of $x\in\RR$, for which there are infinitely
many polynomials $P\in P_n$ satisfying the system
\begin{equation}\label{e:022}
  |P(x)| < H(P)^{-n+1}, \ |P'(x)|< H(P)^{-\delta},
\end{equation}
has zero measure. Indeed, this set consists of points $x\in I_0$
which belong to infinitely many sets
$$
\tau_m=\{x\in I_0:\text{(\ref{e:022}) holds for some }P\in P_n\text{
with }2^{m-1}<H(P)\le 2^m\}
$$
By Theorem~\ref{theorem2}, $|\tau_m|\ll 2^{-m\delta'}$ with
$\delta'>0$. Therefore, $\sum_{m=1}^\infty|\tau_m|<\infty$ and the
Borel-Cantelli Lemma completes the proof of the claim.

Taking into account (\ref{e:006}), this completes the proof of the
case of a small derivative and the proof of
Theorem~\ref{theorem1}.

\mysection{Concluding remarks}
An analogue of Theorem~\ref{theorem1} when $P$ is assumed to be
irreducible over $\QQ$ and primitive (\ie with coprime
coefficients) can also be sought. To make it more precise, let
$P^*_n(H)$ be the subset of $P_n(H)$ consisting of primitive
irreducible polynomials $P$ of degree $\deg P=n$ and height
$H(P)=H$. Now the set of primitive irreducible polynomials of
degree $n$ is $P^*_n=\bigcup_{H=1}^\infty P^*_n(H)$. Let $\wirr$
be the set of $x\in\RR$ such that there are infinitely many $P\in
P^*$ satisfying (\ref{e:001}).

\begin{theorem}\label{theorem3}
Let $\Psi:\RR\to\RR^+$ be arbitrary function such that the sum
\begin{equation}\label{e:023}
\sum_{H=1}^\infty \frac{\# P^*_n(H)}{H}\,\Psi(H)
\end{equation}
converges. Then $|\wirr|=0$.
\end{theorem}

For $n=1$ the proof of Theorem~\ref{theorem3} is a straightforward
application of the Borel-Cantelli Lemma and we again refer to
\cite[p.\,121]{Cassels1}. For $n>1$ the proof follows from the
following 2 observations: 1) $\wirr\subset\w$ and 2) $\#
P^*_n(H)\asymp H^n$. The second one guarantees the converges of
(\ref{e:002}), which now implies $0\le|\wirr|\le|\w|=0$. The proof
of the relation $\# P^*_n(H)\asymp H^n$ is elementary and is left
for the reader. In fact, one can easily estimate the number of
primitive reducible polynomials in $P_n$ and take them off the set
of all primitive polynomials in $P_n$ which is well known to
contain at least constant $\times\#P_n(H)$.

\smallskip

{\textsf{The Duffin-Schaeffer conjecture}}. The conjecture states
that for $n=1$ if (\ref{e:023}) diverges then
$|\RR\setminus\wirr|=0$. The multiple $\#P_1^*(H)$ in sum
(\ref{e:023}) becomes $\asymp\varphi(H)$, where $\varphi$ is the
Euler function.

The following problem can be regarded as the generalization of the
Duffin-Schaeffer conjecture for integral polynomials of higher
degree:
$$
 \mbox{{\sl Prove that $|\RR\setminus\wirr|=0$
whenever (\ref{e:023}) diverges.}}
$$

Alternatively,  for $n>1$ one might investigate the measure of
$\RR\setminus\w$. So far it is unclear if for $n>1$
$|\RR\setminus\wirr|=0$ is equivalent to $|\RR\setminus\w|=0$,
which is another intricate question.

\smallskip

{\textsf{A remark on manifolds.}} In the metric theory of
Diophantine approximation on manifolds one usually studies sets of
$\Psi$-approximable points lying on a manifold with respect to the
measure induced on that manifold. Mahler's problem and its
generalisations can be regarded as Diophantine approximation on
the Veronese curve $(x,x^2,\dots,x^n)$.

A point $\vv f\in\RR^n$ is called\/ {\sl $\Psi$-approximable}\/ if
$$
  \|\vv a\cdot\vv f\|<\Psi(|\vv a|_\infty),
$$
for infinitely many $\vv a\in\ZZ^n$, where $|\vv
a|_\infty=\max_{1\le i\le n}|a_i|$ for $\vv a=(a_1,\dots,a_n)$,
$\|x\|=\min\{|x-z|:z\in\ZZ\}$ and $\Psi:\RR\to\RR^+$.

Let $\vv f:U\to\RR^n$ be a map defined on  an open set
$U\subset\RR^d$. We say that $\vv f$ is {\sl non-degenerate at}\/
$\vv x_0\in U$ if for some $l\in\NN$ the map $\vv f$ is $l$ times
continuously differentiable on a sufficiently small ball centered
at $\vv x_0$ and there are $n$ linearly independent over $\RR$
partial derivatives of $\vv f$ at $\vv x_0$ of orders up to $l$.
We say that $\vv f$ is {\sl non-degenerate}\/ if it is
non-degenerate almost everywhere on $U$. The non-degeneracy of a
manifold is naturally defined via the non-degeneracy of its local
parameterisation.

In 1998 D.~Kleinbock and D.~Margulis proved the Baker-\Sprindzuk
conjecture by showing that any non-degenerate manifold is strongly
extremal. In particular, this implies an analogue of Mahler's
problem for non-degenerate manifolds. A few years later an
analogue of A.~Baker's conjecture with monotonic $\Psi$ (normally
called a Groshev type theorem for convergence) has independently
been proven by V.~Beresnevich \cite{Beresnevich-2002a} and by
V.~Bernik, D.~Kleinbock and G.~Margulis
\cite{BernikKleinbockMargulis-2001a} for non-degenerate manifolds.
It is also remarkable that the proofs were given with different
methods. The divergence counterpart (also for monotonic $\Psi$)
has been established in
\cite{BeresnevichBernikKleinbockMargulis-2002}. In
\cite{BernikKleinbockMargulis-2001a} a multiplicative version of
the Groshev type theorem for convergence has also been given.

Theorem~\ref{theorem1} of this paper can be readily generalised
for non-degenerate curves: {\it Given a non-degenerate map $\vv
f:I\to\RR^n$ defied on an interval $I$, for any function
$\Psi:\RR\to\RR^+$ such that the sum\/ {\rm(\ref{e:002})}\/
converges for almost all $x\in I$ the point $\vv f(x)$ is not
$\Psi$-approximable. } Even further, using the slicing technique
of Pyartly \cite{Pyartli-1969} one can extend this for a class of
$n$-differentiable non-degenerate manifolds which can be foliated
by non-degenerate curves. In particular, this class includes
arbitrary non-degenerate analytic manifold. However with the
technique in our disposal we are currently unable to prove the
following

{\sc Conjecture.}  Let $\vv f:U\to\RR^n$ be a non-degenerate map,
where $U$ is an open subset of $\RR^d$. Then for any function
$\Psi:\RR\to\RR^+$ such that the sum\/ {\rm(\ref{e:002})}\/
converges for almost all $\vv x\in U$ the point $\vv f(\vv x)$ is
not $\Psi$-approximable.



\providecommand{\bysame}{\leavevmode\hbox
to3em{\hrulefill}\thinspace}
\providecommand{\MR}{\relax\ifhmode\unskip\space\fi MR }
\providecommand{\MRhref}[2]{%
  \href{http://www.ams.org/mathscinet-getitem?mr=#1}{#2}
} \providecommand{\href}[2]{#2}

\noindent\it
Institute of Mathematics\\
Academy of Sciences of Belarus\\
220072, Surganova 11, Minsk, Belarus\\
{\tt beresnevich@im.bas-net.by}

\end{document}